\documentclass[a4paper,12pt]{article}
\usepackage{graphicx}
\usepackage{a4wide}
\usepackage[english]{babel}
\usepackage{amsfonts}
\usepackage{amsmath}
\usepackage{amssymb}
\usepackage{dsfont}
\usepackage{amsthm}
\newtheorem{lemma}{Lemma}

\newtheorem{theorem}{Theorem}

\allowdisplaybreaks

\newcommand {\E} {\mathbb{E}}
\DeclareMathOperator {\var}{Var_{[0,1]}}

\newcommand {\N} {\mathbb{N}}

\newcommand {\F} {\mathcal{F}}

\makeatletter\def\blfootnote{\xdef\@thefnmark{}\@footnotetext}
\makeatother

\title{\bf On the asymptotic behavior of weakly lacunary series}

\author{C.\ Aistleitner\footnote{Graz University of Technology, Department for Analysis and Computational Number Theory, Steyrergasse 30, 8010 Graz, Austria \mbox{e-mail}: \texttt{aistleitner@math.tugraz.at}.
Research supported by FWF grant S9603-N23.}, I.\ Berkes\footnote{ Graz University of Technology, Institute of
Statistics, M\"unzgrabenstra{\ss}e 11, 8010 Graz, Austria.  \mbox{e-mail}: \texttt{berkes@tugraz.at}. Research
supported by the FWF Doctoral Program on Discrete Mathematics (FWF DK W1230-N13), FWF grant S9603-N23 and OTKA grants K 67961 and K 81928.} and R.\ Tichy\footnote{Graz University of Technology, Department for Analysis and Computational Number Theory, Steyrergasse 30, 8010 Graz, Austria. \mbox{e-mail}: \texttt{tichy@tugraz.at}.
Research supported by the FWF Doctoral Program on Discrete Mathematics (FWF DK W1230-N13) and FWF grant S9603-N23.}}

\begin{document}


\date{}
\maketitle

\abstract{\medskip Let $f$ be a measurable function satisfying
$$f(x+1)=f(x), \quad \int_0^1 f(x)\, dx=0, \quad \var f < + \infty,
$$
and let $(n_k)_{k\ge 1}$ be a sequence of integers satisfying
$n_{k+1}/n_k \ge q >1$ $(k=1, 2, \ldots)$. By the classical theory
of lacunary series, under suitable Diophantine conditions on
$n_k$, $(f(n_kx))_{k\ge 1}$ satisfies the central limit theorem
and the law of the iterated logarithm. These results extend for a
class of subexponentially growing sequences $(n_k)_{k\ge 1}$ as
well, but as Fukuyama (2009) showed, the behavior of $f(n_kx)$ is
generally not permutation-invariant, e.g.\ a rearrangement of the
sequence can ruin the CLT and LIL. In this paper we construct an
infinite order Diophantine condition implying the
permutation-invariant CLT and LIL without any growth conditions on
$(n_k)_{k\ge 1}$ and show that the known finite order Diophantine
conditions in the theory do not imply  permutation-invariance even
if $f(x)=\sin 2\pi x$ and $(n_k)_{k\ge 1}$ grows almost
exponentially. Finally we prove that, in a suitable statistical
sense, for almost all sequences $(n_k)_{k\ge 1}$ growing faster than polynomially, $(f(n_kx))_{k\ge 1}$ has  permutation-invariant behavior.}

\section{Introduction}

Let $f$ be a measurable function satisfying
\begin{equation}\label{fcond}
f(x+1)=f(x), \quad \int_0^1 f(x)\, dx=0, \quad \var f < + \infty
\end{equation}
and let $(n_k)_{k\ge 1}$ be a sequence of positive integers
satisfying the Hadamard gap condition
\begin{equation}\label{had}
n_{k+1}/n_k\ge q>1 \qquad(k=1, 2, \ldots).
\end{equation}
In the case $n_k=2^k$, Kac \cite{kac46} proved that  $f(n_kx)$
satisfies the central limit theorem
\begin{equation}\label{fclt}
N^{-1/2}\sum_{k=1}^N f(n_k x) \stackrel{\mathcal D}
{\longrightarrow} \mathcal{N}(0, \sigma^2)\\
\end{equation}
with respect to the probability space $[0, 1]$ equipped with the
Lebesgue measure, where
$$
\sigma^2=\int_0^1 f^2(x)\, dx +2\sum_{k=1}^\infty \int_0^1
f(x)f(2^kx)\, dx.
$$
Gaposhkin \cite{gap66} extended  (\ref{fclt}) to  the case when
the fractions $n_{k+1}/n_k$ are all integers or if
$n_{k+1}/n_k\to\alpha$, where $\alpha^r$ is irrational for $r=1,
2, \ldots$. On the other hand, an example of Erd\H{o}s and Fortet
(see \cite{kac49}, p.\ 646) shows that the CLT (\ref{fclt}) fails
if $n_k=2^k-1$. Gaposhkin also showed (see \cite{gap70}) that the
asymptotic behavior of $\sum_{k=1}^N f(n_kx)$ is intimately
connected with the number of solutions of the Diophantine equation
$$
an_k+bn_l=c, \qquad 1\le k, l\le N.
$$
Improving these results, Aistleitner and Berkes \cite{aibe} gave a
necessary and sufficient condition for the CLT (\ref{fclt}). For
related laws of the iterated logarithm, see \cite{beph},
\cite{gap66}, \cite{iz}, \cite{my}.

The previous results show that for arithmetically "nice" sequences
$(n_k)_{k\ge 1}$, the system $f(n_kx)$ behaves like a sequence of
independent random variables. However, as an example of Fukuyama
\cite{ft2} shows,  this result is not permutation-invariant: a
rearrangement of $(n_k)_{k\ge 1}$ can change the variance of the
limiting Gaussian law or ruin the CLT altogether. A complete
characterization of the permutation-invariant CLT and LIL for
$f(n_kx)$ under the Hadamard gap condition (\ref{had}) is given in
our forthcoming paper \cite{aibeti}. In particular, it is shown
there that in the harmonic case $f(x)=\cos 2\pi x$, $f(x)=\sin
2\pi x$ the CLT and LIL for $f(n_kx)$ hold after any permutation
of $(n_k)_{k \ge 1}$.

For subexponentially growing $(n_k)_{k\ge 1} $ the situation
changes radically. Note that in the case $f(x)=\cos 2\pi x$,
$f(x)=\sin 2\pi x$, the  unpermuted CLT and LIL remain valid under
the weaker gap condition
$$ n_{k+1}/n_k\ge 1+ck^{-\alpha}, \qquad 0<\alpha<1/2,
$$
see Erd\H{o}s \cite{er62}, Takahashi \cite{tak65}, \cite{tak72}.
However, as the following theorem shows, the slightest weakening
of the Hadamard gap condition (\ref{had}) can ruin the
permutation-invariant CLT and LIL.

\begin{theorem}\label{counterex} For any positive sequence
$(\varepsilon_k)_{k \geq 1}$ tending to 0, there exists a sequence
$(n_k)_{k \geq 1}$ of positive integers satisfying
\begin{equation}\label{epsilongap}
n_{k+1}/n_k\ge 1+\varepsilon_k, \qquad k\ge k_0
\end{equation}
and a permutation $\sigma: {\mathbb N}\to {\mathbb N}$ of the
positive integers such that
\begin{equation}\label{ell1}
N^{-1/2} \sum_{k=1}^N \cos 2\pi n_{\sigma(k)} x -b_N
\stackrel{\mathcal D}
\longrightarrow G\\
\end{equation}
where $G$ is a nongaussian distribution with characteristic
function given by (\ref{infdiv})-(\ref{FG}) and $(b_N)_{N\ge 1}$
is a numerical sequence with $b_N=O(1)$. Moreover, there exists a
permutation $\sigma: {\mathbb N}\to {\mathbb N}$ of the positive
integers such that
\begin{equation} \label{lilcount}
\limsup_{N \to \infty} \frac{\sum_{k=1}^N \cos 2 \pi n_{\sigma(k)}
x} {\sqrt{2 N \log \log N}} = + \infty \quad \textup{a.e.}
\end{equation}
\end{theorem}

Even the number theoretic conditions implying the CLT and LIL
under subexponential gap conditions do not help here: the sequence
$(n_k)_{k\ge 1}$ in Theorem \ref{counterex} can be chosen so that
it satisfies conditions {\bf B}, {\bf C}, {\bf G} in our paper
\cite{bpt} implying  very strong independence properties of $\cos
2\pi n_kx$,  $\sin 2\pi n_kx$, including the CLT and LIL. In fact,
it is very difficult to construct subexponential sequences
$(n_k)_{k\ge 1}$ satisfying the permutation-invariant CLT and LIL:
the only known example (see \cite{HLP}) is the
Hardy-Littlewood-P\'olya sequence, i.e.\ the sequence generated by
finitely many primes and arranged in increasing order; the proof
uses deep number theoretic tools. The purpose of this paper is to
introduce a new, infinite order Diophantine condition
$\bf{A}_{\boldsymbol\omega}$ which implies the
permutation-invariant CLT and LIL for $f(n_kx)$ and then to show that,
in a suitable statistical sense, almost all sequences $(n_k)_{\ge 1}$ growing faster than polynomially satisfy $\boldsymbol{A_\omega}$.
Thus, despite the difficulties to construct explicit examples, the
permutation-invariant CLT and LIL are rather the rule than the
exception.

Given a nondecreasing sequence $\boldsymbol{\omega}=(\omega_1,
\omega_2, \ldots)$ of positive numbers  tending to $+\infty$, let
us say that a sequence $(n_k)_{k\ge 1}$ of different positive
integers satisfies

\bigskip\noindent
{\bf Condition $\boldsymbol{A_\omega}$}, if for any $N\ge N_0$ the
Diophantine equation
\begin{equation}\label{aomegadef}
a_1 n_{k_1} +\ldots +a_r n_{k_r}=0, \qquad \ 2\le r\le \omega_N,
 \ 0< |a_1|, \ldots, |a_r| \le  N^{\omega_N}
\end{equation}
with different indices $k_j$ and nonzero integer coefficients
$a_j$ has only such solutions where all  $n_{k_j}$ belong to the smallest
$N$ elements of the sequence $(n_k)_{k\ge 1}$. \\

Clearly, this property is permutation-invariant and it implies that for
any fixed nonzero integer coefficients $a_j$ the number of solutions of
(\ref{aomegadef}) with different indices $k_j$ is at most $N^r$.

\begin{theorem}\label{th1} Let {\boldmath{$\omega$}}=$(\omega_1,
\omega_2, \ldots)$ be a nondecreasing sequence tending to
$+\infty$ and let $(n_k)_{k\ge 1}$ be a sequence of different
positive integers satisfying condition $\boldsymbol{A_\omega}$.
Then for any $f$ satisfying (\ref{fcond}) we have
\begin{equation}\label{permclt}
N^{-1/2}\sum_{k=1}^N f(n_k x) \stackrel{\mathcal D}
{\rightarrow} \mathcal{N}(0, \|f\|^2)
\end{equation}
where  $\|f\|$ denotes the $L_2(0, 1)$ norm of $f$.
If $\omega_k \ge (\log k)^\alpha$ for some $\alpha>0$
and $k\ge k_0$, then we also have
\begin{equation} \label{permlil}
\limsup_{N \to \infty} \frac{\sum_{k=1}^N f(n_k x)}{(2 N
\log \log N)^{1/2}} = \|f\| \quad \textup{a.e.}
\end{equation}
\end{theorem}

Condition $\boldsymbol{A_\omega}$ is different from the usual
Diophantine conditions in lacunarity theory, which typically
involve 4 or less terms.
In contrast, $\boldsymbol{A_\omega}$ is an 'infinite order'
condition, namely it involves equations with arbitrary large
order. As noted, the usual Diophantine conditions do not suffice
in Theorem \ref{th1}. Given any $\omega_k\uparrow\infty$, it is
not hard to see that any sufficiently rapidly growing sequence
$(n_k)_{k\ge 1}$ satisfies $\boldsymbol{A_\omega}$; on the other
hand, we do not have any "concrete" subexponential examples for
$\boldsymbol{A_\omega}$. However, we will show that, in a suitable
statistical sense, almost all sequences growing faster than
polynomially satisfy condition $\boldsymbol{A_\omega}$ for some
appropriate $\boldsymbol{\omega}$. To make this precise requires
defining a probability measure over the set of such sequences, or,
equivalently, a natural random procedure to generate such
sequences. A simple procedure is to choose $n_k$ independently and
uniformly from the integers in the interval
\begin{equation}\label{I}
I_k=[a(k-1)^{\omega_{k-1}}, a k^{\omega_k}), \qquad k=1, 2, \ldots.
\end{equation}
Note that the length of $I_k$ is at least $a\omega_k
(k-1)^{\omega_k-1}\ge a\omega_1$ for $k=2, 3, \ldots $ and equals
$a$ for $k=1$ and thus choosing $a$ large enough, each $I_k$
contains at least one integer. Let $\mu_{\boldsymbol \omega}$ be
the distribution of the random sequence $(n_k)_{k \ge 1}$ in the
product space $I_1\times I_2 \times \ldots$.

\begin{theorem} \label{th3} Let $\omega_k\uparrow \infty$ and let
$f$ be a function satisfying (\ref{fcond}). Then with probability
one with respect to $\mu_{\boldsymbol \omega}$ the sequence
$(f(n_k x))_{k \geq 1}$ satisfies the CLT (\ref{permclt}) after
any permutation of its terms, and if $\omega_k\ge (\log
k)^\alpha$ for some $\alpha>0$ and $k\ge k_0$, $(f(n_k x))_{k \geq 1}$
also satisfies the LIL (\ref{permlil}) after any permutation of
its terms.\\
\end{theorem}

The sequences $(n_k)_{k\ge 1}$ provided by $\mu_{\boldsymbol
\omega}$ satisfy $n_k=O(k^{\omega_k})$; for slowly increasing
$\omega_k$ the so obtained sequences grow much slower than
exponentially, in fact they grow barely faster than polynomial
speed. If $\omega_k$ grows so slowly that
$\omega_k-\omega_{k-1}=o((\log k)^{-1})$, then the so obtained
sequence $(n_k)_{k\ge 1}$ has the
precise speed $n_k \sim k^{\omega_k}$. We do not know if there exist
polynomially growing sequences $(n_k)_{k \geq 1}$ satisfying the
permutation-invariant CLT or LIL. The proof of Theorem 3 will also
show that with probability 1, the sequences provided by
$\mu_{\boldsymbol \omega}$ satisfy $\boldsymbol{A_{\omega^*}}$ with
${\boldsymbol \omega^*}= (c\omega_1^{1/2}, c\omega_2^{1/2}, \ldots)$.

\section{Proofs}

\subsection{Proof of Theorem \ref{counterex}}

We begin with the CLT part. Let $(\varepsilon_k)_{k\ge 1}$ be a
positive sequence tending to 0. Let $m_1<m_2< \ldots$ be positive
integers such that $m_{k+1}/m_k\ge 2^{k^2}$, $k=1, 2, \ldots$ and
all the $m_k$ are powers of 2; let $r_1\le r_2\le \ldots$ be
positive integers satisfying $1\le r_k\le k^2$. Put $I_k=\{m_k,
2m_k, \ldots, r_km_k\}$; clearly the sets $I_k$, $k=1,2,...$ are
disjoint. Define the sequence $(n_k)_{k\ge 1}$ by
\begin{equation}\label{nk}
(n_k)_{k\ge 1}=\bigcup\limits_{j=1}^\infty I_j.
\end{equation}
Clearly, if $n_k, n_{k+1}\in I_j$, then $n_{k+1}/n_k\ge 1+1/r_j$
and thus if $r_j$ grows sufficiently slowly, the sequence
$(n_k)_{k\ge 1}$ satisfies the gap condition (\ref{epsilongap}).
Also, if $r_j$ grows sufficiently slowly, there exists a
subsequence $(n_{k_\ell})_{\ell\ge 1}$ of $(n_k)_{k\ge 1}$ which
has exactly the same structure as the sequence in (\ref{nk}), just
with $r_k\sim k$. By the proof of Theorem~1 in \cite{ber91},
$(\cos 2 \pi n_{k_\ell} x)_{\ell \geq 1}$ satisfies
\begin{equation}\label{cauchy}
\frac{1}{\sqrt{N}} \sum_{\ell=1}^N \cos 2\pi n_{k_\ell} x- b_N
\stackrel{\mathcal D} {\rightarrow} G
\end{equation}
where $(b_N)_{N\ge 1}$ is a numerical sequence with $b_N=O(1)$ and
$G$ is a nongaussian infinitely divisible distribution with
characteristic function
\begin{equation}\label{infdiv}
\exp\left\{\int\limits_{R\backslash\{0\}}\left(e^{itx}-1-
{itx\over{1+x^2}}\right)dL(x)\right\}
\end{equation}
where
\begin{equation}\label{L} \qquad L(x)= \begin{cases}
-{1\over{\pi}}\int\limits^1_{x}{F(t)\over{t}}dt&  \hbox{if} \qquad 0<x\le
1 \cr {1\over{\pi}}\int\limits^1_{-x}{G(t)\over{t}}dt & \hbox{if} \qquad
-1\le x<0 \cr
       0 & \hbox{if} \qquad  \mid x\mid >1
       \end{cases}
\end{equation}
and \begin{equation}\label{FG}
 F(t)=\lambda \{x>0:\sin x/x\ge t\},
\quad G(t)=\lambda \{ x>0:\sin x/x\le -t\} \quad (t>0)
\end{equation}
where $\lambda$ is the Lebesgue measure. Define a permutation
$\sigma$ in the following way:
\begin{itemize}
\item for $k\notin \{1, 2, 4, \ldots, 2^m, \ldots\}$ \  $\sigma(k)$ takes
the values of the set $\{k_1, k_2, \ldots \}$ in consecutive order
\item for $k\in \{1, 2, 4, \ldots, 2^m, \ldots\}$ \
$\sigma(k)$ takes the values of the set ${\mathbb N} \setminus
\{k_1, k_2, \ldots\}$ in consecutive order.
\end{itemize}
Then $\sigma$ is a permutation of $\mathbb N$ and the sums
$$
\sum_{k=1}^N \cos 2 \pi n_{\sigma(k)} x \qquad \textrm{and} \qquad
\sum_{l=1}^N \cos 2 \pi n_{k_l} x
$$
differ at most in $2 \log_2 N$ terms. Therefore, (\ref{cauchy})
implies (\ref{ell1}), proving the first part Theorem 1.\\

The proof of the LIL part of Theorem 1 is modeled after the proof
of Theorem 1 in Berkes and Philipp \cite{bp94}. Similarly as
above, we construct a sequence $(n_k)_{k \geq 1}$ satisfying
(\ref{epsilongap}) that contains a subsequence $(\mu_k)_{k\ge 1}$
of the form (\ref{nk}) with $I_k=\{m_k, 2m_k, \ldots, r_k m_k\}$,
where $r_k \sim k \log k$ and $(m_k)_{k \geq 1}$ is growing fast;
specifically we choose $m_k$ in such a way that it is a power of
$2$ and $m_{k+1} \geq r_k 2^{2k} m_{k}$.
Let $\F_i$ denote the $\sigma$-field generated of the dyadic
intervals
$$
[\nu 2^{-(\log_2 m_i)-i},(\nu+1)2^{-(\log_2 m_i)-i}), \quad 0 \leq
\nu < 2^{(\log_2 m_i)+i}.
$$
Write
$$
X_i = \cos 2 \pi m_i x + \dots + \cos 2 \pi r_i m_i x
$$
and
$$
Z_i = \E(X_i|\F_i).
$$
Then for all $x \in (0,1)$
$$
|X_i(x) - Z_i(x)| \ll r_i^2 m_i 2^{-(\log_2 m_i)-i}
$$
and thus
$$
\sum_{i \geq 1} |X_i(x)-Z_i(x)| < \infty \qquad \textrm{for all} \
x\in (0,1).
$$
Like in \cite[Lemma 2.1]{bp94} the random variables $Z_1, Z_2,
\dots$ are independent, and like in \cite[Lemma 2.2]{bp94} for
almost $x \in (0,1)$ we have
$$
\limsup_{i \to \infty} \, X_i/m_i \geq 2/\pi.
$$
Assume that for a fixed $x$ and some $i\ge 1$ we have $X_i/m_i
\geq 1/\pi$. Then either
$$
|X_1 + \dots + X_{i-1}| \geq m_i/2 \pi,
$$
or
$$
|X_1+ \dots+X_i| \geq m_i/2\pi.
$$
Since the total number of summands in $X_1, \dots,X_i$ is $\ll i^2
\log i$, and since
$$m_i \gg (i^2 \log i)^{1/2} (\log(i^2 \log i))^{1/2},$$
we have
$$
\limsup_{N \to \infty} \frac{\left| \sum_{k=1}^N \cos 2 \pi \mu_k
x \right|}{\sqrt{N \log N}} > 0 \qquad \textup{a.e.},
$$
and, in particular,
$$
\limsup_{N \to \infty} \frac{\left| \sum_{k=1}^N \cos 2 \pi \mu_k
x \right|}{\sqrt{2 N \log \log N}} = + \infty \qquad
\textup{a.e.},
$$
Thus we constructed a subsequence of $(n_k)_{k\ge 1}$ failing the
LIL and similarly as above, we can construct a permutation
$(n_{\sigma(k)})_{k\ge 1}$ of $(n_k)$ failing the LIL as well.

\subsection{Proof of Theorem \ref{th1}}

\begin{lemma} \label{lemma1}
Let $\omega_k\uparrow \infty$ and let $(n_k)_{k \geq 1}$ be a sequence
of different positive integers satisfying condition
$\boldsymbol{A_\omega}$. Let $f$ satisfy (\ref{fcond}) and  put
$S_N=\sum_{k=1}^N f(n_k x)$, $\sigma_N=(\E S_N^2)^{1/2}$. Then for
any $p\ge 3$ we have
$$
\E S_N^p=
\begin{cases}
\frac{p!}{(p/2)!}2^{-p/2} \sigma_N^p+O(T_N)
 \quad &\text{if} \ p \ \text{is even} \\
O(T_N) \quad &\text{if} \ p \ \text{is odd}
\end{cases}
$$
where
\begin{equation*}
T_N= \exp (p^2) N^{(p-1)/2}(\log N)^p
\end{equation*}
and the constants implied by the $O$ are absolute.
\end{lemma}

\noindent {\bf Proof.} Fix $p\ge 2$ and choose the integer $N$ so
large that $\omega_{[N^{1/4}]}\ge 8p$. Without loss of generality
we may assume that $f$ is an even function and that $\|f\|_\infty
\le 1$, $\var f \leq 1$; the proof in the general case is similar. Let
\begin{equation} \label{frep}
f \sim \sum_{j=1}^\infty a_j \cos 2 \pi j x
\end{equation}
be the Fourier series of $f$. $\var f \leq 1$ implies
\begin{equation} \label{fouriercoeffs}
|a_j| \leq j^{-1},
\end{equation}
(see Zygmund \cite[p.\ 48]{zt}) and writing
$$
g(x) = \sum_{j=1}^{N^{p}} a_j \cos 2 \pi j x, \qquad r(x) = f(x) -
g(x),
$$
we have
$$
\|g\|_\infty \leq \var f + \|f\|_\infty \leq 2, \qquad
\|r\|_\infty \leq  \|f\|_\infty + \|g\|_\infty \leq 3
$$
by (4.12) of Chapter II and (1.25) and (3.5) of Chapter III of
Zygmund \cite{zt}. Letting $\|\cdot\|$ and $\| \cdot \|_p$ denote
the $L_2(0,1)$, resp. $L_p(0,1)$ norms, respectively,
(\ref{fouriercoeffs}) yields for any positive integer $n$
\begin{equation} \label{likein2}
\|(r(n x)\|^2 = \|r(x)\|^2 = \frac{1}{2} \sum_{j=N^{p}+1}^\infty
a_j^2 \leq N^{-p}.
\end{equation}
By Minkowski's inequality,
\begin{equation*}\label{min1}
\|S_N\|_p \leq \| \sum_{k=1}^N g(n_k x) \|_p +  \| \sum_{k=1}^N
r(n_k x) \|_p,
\end{equation*}
and
\begin{equation}\label{min2}
\| \sum_{k=1}^N r(n_k x) \|_p  \leq 3\sum_{k=1}^N \left\| r(n_k
x)/3 \right\|_p  \leq 3 \sum_{k=1}^N \left\| r(n_k x)/3
\right\|^{2/p} \leq 3 \sum_{k=1}^N N^{-1} \leq 3.
\end{equation}
Similarly,
$$
\left| \|S_N\| - \| \sum_{k=1}^N g(n_k x) \| \right| \leq \|
\sum_{k=1}^N r(n_k x) \| \leq N^{-\frac{p}{2}+1}
$$
and therefore
\begin{eqnarray}
& & \left| \|S_N\|^p - \| \sum_{k=1}^N g(n_k x) \|^p
\right| \nonumber\\
& \leq & p ~\max\left(\|S_N\|^{p-1}, \|\sum_{k=1}^N g(n_k x)
\|^{p-1} \right) \cdot \left|\|S_N\| - \|
\sum_{k=1}^N g(n_k x) \|\right| \nonumber\\
& \ll & p \left(N (\log \log N)^2 \right)^{\frac{p-1}{2}}
N^{-\frac{p}{2}+1} \nonumber\\
& \ll & p (\log \log N)^{p-1} N^{1/2} \label{diffsg}
\end{eqnarray}
since by a result of G\'al \cite{ga} and Koksma \cite{koa}
\begin{equation} \label{koks}
\|S_N\|^2 \ll N (\log \log N)^2 \quad \textrm{and} \quad
\|\sum_{k=1}^N g(n_k x) \|^2 \ll N (\log \log N)^2,
\end{equation}
where the implied constants are absolute.\\

By expanding and using elementary properties of the trigonometric
functions we get
\begin{align}\label{esnp}
&\E\left(\sum_{k=1}^N g(n_k x)\right)^p\nonumber\\
&= 2^{-p} \sum_{1\le j_1, \ldots, j_p \le N^{p}} a_{j_1}\cdots
a_{j_p}\sum_{1\le k_1, \ldots, k_p\le N} ~\mathds{I} \{\pm j_1
n_{k_1} \pm \ldots \pm j_p n_{k_p}=0\},
\end{align}
with all possibilities of the signs $\pm$ within the indicator
function.  Assume that $j_1, \dots, j_p$ and the signs $\pm$ are
fixed, and consider a solution of  $\pm j_1 n_{k_1} \pm \ldots \pm
j_p n_{k_p}=0$.  Then the set $\{1, 2, \ldots, p\}$ can be split
into disjoint sets $A_1, \ldots, A_l$ such that for each such set
$A$ we have $\sum_{i\in A} \pm j_i n_{k_i}=0$ and no further
subsums of these sums are equal to 0. Group the terms of
$\sum_{i\in A} \pm j_i n_{k_i}$ with equal $k_i$. If after
grouping there are at least two terms, then by the restriction on
subsums, the sum of the coefficients $j_i$ in each group will be
different from 0 and will not exceed
$$
pN^{p} \le \omega_{[N^{1/4}]}N^{\frac{1}{8} \omega_{[N^{1/4}]}}
\le  2^{\frac{1}{8}\omega_{[N^{1/4}]}} N^{\frac{1}{8}
\omega_{[N^{1/4}]}}\le N^{\frac{1}{4} \omega_{[N^{1/4}]}} \qquad
(N\ge N_0).
$$
Also the number of terms after grouping will be at most $p\le
\omega_{[N^{1/4}]}$   and thus applying condition
$\boldsymbol{A_\omega}$ with the index $[N^{1/4}]$ shows that
within a block $A$ the
$n_{k_i}$ belong to the smallest $[N^{1/4}]$ terms of the
sequence. Thus letting $|A|=m$, the number of solutions of
$\sum_{i\in A} \pm j_i n_{k_i}=0$ is at most $N^{m/4}$. If after
grouping there is only one term, then all the $k_i$ are equal and
thus the number of solutions of $\sum_{i\in A} \pm j_i n_{k_i}=0$
is at most $N$. Thus if $m\ge 3$ then the number of solutions of
$\sum_{i\in A} \pm j_i n_{k_i}=0$ in a block is at most $N^{m/3}$.
If $m=2$, then the number of solutions is clearly at most $N$.
Thus if $s_i=|A_i|$ \ ($1\le i \le l$) denotes the cardinality of
$A_i$, the number of solutions of $\pm j_1 n_{k_1} \pm \ldots \pm
j_p n_{k_p}=0$ admitting such a decomposition with fixed $A_1,
\ldots, A_l$ is at most
\begin{align*}
&\prod_{\{i: s_i\ge 3\}} N^{s_i/3}  \prod_{\{i: s_i=2\}} N =
N^{\frac{1}{3} \sum_{\{i: s_i\ge 3\}}s_i+ {\sum_{\{i: s_i=2\}}1}} \\
& = N^{\frac{1}{3}\sum_{\{i: s_i\ge 3\}}s_i+
\frac{1}{2}{\sum_{\{i: s_i=2\}}s_i}}
= N^{\frac{1}{3}\sum_{\{i: s_i\ge 3\}}s_i +\frac{1}{2}(p-
\sum_{\{i: s_i\ge
3\}}s_i)}\\
&=N^{\frac{p}{2}-\frac{1}{6}\sum_{\{i: s_i\ge 3\}}s_i}.
\end{align*}
If there is at least one $i$ with $s_i\ge 3$, then the last
exponent is at most $(p-1)/2$ and since the number of partitions
of the set $\{1, \ldots, p\}$ into disjoint subsets is at most
$p!\, 2^p$, we see that the number of solutions of $\pm j_1
n_{k_1} \pm \ldots \pm j_p n_{k_p}=0$ where at least one of the
sets $A_i$ has cardinality $\ge 3$ is at most $p!\, 2^p
N^{(p-1)/2}$. If $p$ is odd, there are no other solutions and thus
using (\ref{fouriercoeffs}) the inner sum in (\ref{esnp}) is at
most $p!\, 2^p N^{(p-1)/2}$ and consequently, taking into account
the $2^p$ choices for the signs $\pm 1$,
\begin{eqnarray*}\label{esnpodd}
& & \left|\E\left(\sum_{k\le N} g(n_kx)\right)^p\right| \nonumber\\
&\le & p!\, 2^pN^{(p-1)/2} 2^p \sum_{1\le j_1, \ldots, j_p \le
N^{p}} |a_{j_1}\cdots a_{j_p}| \ll \exp(p^2) N^{(p-1)/2} (\log
N)^p.
\end{eqnarray*}
If $p$ is even, there are also solutions where each $A$ has
cardinality 2. Clearly, the contribution of the terms in
(\ref{esnp}) where $A_1=\{1, 2\}, A_2=\{3, 4\}, \ldots$ is
\begin{eqnarray*}
& & \left( \frac{1}{4}\sum_{1\le i,j \le N^{2p}} \sum_{1\le k,
\ell \le N} a_i a_j \mathds{I}\{ \pm i n_k \pm j
n_\ell=0\}\right)^{p/2}  =
\left( \E \left(\sum_{k\le N} g(n_k x) \right)^2 \right)^{p/2} \nonumber \\
& = & \left\| \sum_{k\le N} g(n_kx)\right\|^{p} \nonumber\\
& = & \|S_N\|^p + \mathcal{O} \left( p (\log \log N)^{p-1}
N^{1/2} \right)
\end{eqnarray*}
by (\ref{diffsg}).\\

Since the splitting of $\{1, 2, \ldots, p\}$ into pairs can be
done in $\frac{p!}{(p/2)!} 2^{-p/2}$ different ways, we proved
that
\begin{equation}\label{egnp}
\E \left(\sum_{k\le N} g(n_kx)\right)^p=
\begin{cases}
\frac{p!}{(p/2)!}2^{-\frac{p}{2}} \sigma_N^p
+O(T_N) \\
O (T_N)
\end{cases}
\end{equation}
according as $p$ is even or odd; here
$$
T_N= \exp (p^2) N^{(p-1)/2}(\log N)^p.
$$
Now, letting $G_N=\sum_{k\le N} g(n_kx)$ we get, using
(\ref{min2}), (\ref{koks}) and (\ref{egnp}),
\begin{eqnarray*}\label{sngn}
& &|\E S_N^p-\E G_N^p| \\
& \leq & p ~\max \left( \|S_N\|_p^{p-1}, \|G_N\|_p^{p-1} \right)
\cdot \left| \|S_N\|_p - \|G_N\|_p\right| \\
& \ll & p \left(\frac{p!}{(p/2)!}2^{-\frac{p}{2}}
\sigma_N^p\right)^{\frac{p-1}{p}} \\
& \ll & T_N,
\end{eqnarray*}
completing the proof of Lemma \ref{lemma1}. \qquad $\square$\\

\begin{lemma}\label{sigma}
Let $\omega_k \uparrow \infty$ and let  $(n_k)_{k\ge 1}$ be a
sequence of different positive integers satisfying
condition $\boldsymbol{A_\omega}$. Then for any $f$ satisfying (\ref{fcond})
we have
\begin{equation}\label{sigmacond}
\int_0^1 \left(\sum_{k=1}^N f(n_k x)\right)^2\, dx \sim
\|f\|^2 N \quad \text{as} \ N\to\infty.
\end{equation}
\end{lemma}

\noindent {\bf Proof.} Clearly, $\omega_{[N^{1/4}]} \geq 4$ for
sufficiently large $N$ and thus  applying Condition $\boldsymbol{A_\omega}$
for the index $[N^{1/4}]$ it follows that for $N\ge N_0$ the
Diophantine equation
\begin{equation} \label{nosol1}
j_1 n_{i_1} + j_2 n_{i_2} = 0, \qquad i_1 \neq i_2,\
0< |j_1|,|j_2| \leq N
\end{equation}
has only such solutions where $n_{i_1}, n_{i_2}$ belong to the set $J_N$
of $[N^{1/4}]$ smallest elements  of the sequence $(n_k)_{k\ge 1}$.
Write $p_N(x)$ for the $N$-th partial sum of the Fourier series of $f$,
and $r_N$ for the $N$-th remainder term. Then
we have for any for any $f$ satisfying (\ref{fcond})
\begin{equation}\label{32}
\left\|\sum_{k=1}^N f(n_k x) \right\|  \geq
\left\|\sum_{k\in [1, N]\setminus J_N} p_N(n_k x) \right\| -
\left\|\sum_{k\in J_N} f(n_k x) \right\| -
\left\|\sum_{k\in [1, N] \setminus J_N} r_N(n_k x) \right\|.
\end{equation}
Using the previous remark on the number of solutions of (\ref{nosol1})
we get, as in (\ref{esnp}),
\begin{equation*}
\left\|\sum_{k\in [1, N] \setminus J_N} p_N(n_k x) \right\| =
 (N-[N^{1/4}])^{1/2} ~\|p_N\| \sim N^{1/2} \|f\|,
\end{equation*}
since $\|p_N\|\to \|f\|$. Further, $\|r_N\|\ll N^{-1/2}$ by
(\ref{fouriercoeffs}) and thus using
Minkowski's inequality and the results of G\'al and Koksma
mentioned in (\ref{koks}),  we get
\begin{equation*}
\left\|\sum_{k\in J_N} f(n_k x) \right\|  \ll
N^{1/4}, \quad
\left\|\sum_{k\in [1, N] \setminus J_N} r_N(n_k x) \right\|  \ll
\sqrt{N}\log\log N \|r_N\| \ll \log\log N.
\end{equation*}
These estimates, together with (\ref{32}), prove Lemma \ref{sigma}.
\\

Lemma \ref{lemma1} and Lemma \ref{sigma} imply that for any fixed
$p\ge 2$, the $p$-th moment of $S_N/\sigma_N$ converges to
$\frac{p!}{(p/2)!}2^{-p/2}$ if $p$ is even and to $0$ if $p$ is odd;
in other words, the moments of $S_N/\sigma_N$ converge to the moments
of the standard normal distribution. By $\sigma_N\sim \|f\| \sqrt{N}$
and a well known result in probability theory, this proves the CLT part
of Theorem \ref{th1}. The
proof of the LIL part of Theorem \ref{th1} is more involved, and
we will give just a sketch of the proof. The details can be
modeled after the proof of \cite[Theorem 1]{HLP}. The crucial
ingredient is Lemma \ref{lemmalil} below, which yields the LIL
part of Theorem \ref{th1}, like \cite[Theorem 1]{HLP} follows
from \cite[Lemma 3]{HLP}.\\

Let $\theta>1$ and define $\Delta_M' = \{k \in \N:~\theta^M < k
\leq \theta^{M+1}\}$ and $T_M'=\sum_{k\in \Delta_M'} f(n_kx)$. By
the standard method of proof of the LIL, we need precise bounds
for the tails of $T_M'$ and also, a near independence relation for
the $T_M'$ for the application of the Borel-Cantelli lemma in the
lower half of the LIL. From the set $\Delta_M'$ we remove its
$[\theta^{M/4}]$ elements with the smallest value of $n_k$ (recall
that the sequence $(n_k)_{k\ge 1}$ is not assumed to be
increasing) and denote the remaining set by $\Delta_M$. Since the
number of removed elements is $\ll |\Delta_M'|^{1/4}$, this
operation does not influence the partial sum asymptotics of
$T_M'$. Like in the proof of Lemma \ref{lemma1}, we assume that we
have a representation of $f$ in the form (\ref{frep}) and that
(\ref{fouriercoeffs}) holds. Define
$$
g_M(x)=\sum_{j=1}^{[\theta^M]^2} a_j \cos 2 \pi j x, \qquad
\sigma_M^2 = \int_0^1 \left( \sum_{k \in \Delta_M} g_M(n_k x)
\right)^2 dx
$$
and
$$
T_M= \sum_{k \in \Delta_M} g_M(n_k x), \qquad  Z_M = T_M/\sigma_M.
$$
From Lemma \ref{sigma} it follows easily
\begin{equation}\label{sigmabound}
\sigma_M \gg |\Delta_M|^{1/2}.
\end{equation}
Assume that $(n_k)_{k\ge 1}$ satisfies Condition
$\boldsymbol{A_\omega}$ for a sequence $(\omega_k)_{k\ge 1}$ with
$\omega_k \geq (\log k)^\alpha$ for some $\alpha>0$, $k\ge k_0$.
Without loss of generality we may assume $0< \alpha<1/2$. Choose
$\delta>0$ so small that for sufficiently large $r$
\begin{equation}\label{side}
\left(\log \theta^{\sqrt{r}/4}\right)^\alpha > 4 \left(\log
\theta^{r}\right)^\delta.
\end{equation}
\begin{lemma} \label{lemmalil}
For sufficiently large $M,N$ satisfying $N^{1-\alpha/2} \leq M
\leq N$, and for positive integers $p,q$ satisfying $p+q \leq
(\log \theta^N)^\delta$ we have
$$
\E Z_M^p Z_N^q = \left\{ \begin{array}{ll}
\frac{p!}{(p/2)!2^{p/2}} \frac{q!}{(q/2)!2^{q/2}} + \mathcal{O}
(R_{M,N}) & \textrm{if $p,q$ are even} \\ \mathcal{O} (R_{M,N}) &
\textrm{otherwise}
\end{array} \right.
$$
where
$$
R_{M,N} = 2^{p+q} (p+q)! \, (\log M)^{p+q} |\Delta_M|^{-1/2}.
$$
\end{lemma}

{\bf Proof.} Note that $N\le M^{1+3\alpha/4}$ and thus,  setting
$L=[\theta^{M/4}]$, relation $p+q \leq  (\log \theta^N)^\delta$
and (\ref{side}) imply
\begin{equation*} \label{pq}
p+q \le \frac{1}{4} (\log \theta^{\sqrt{N}/4})^\alpha \le \frac{1}{4} (\log
\theta^{M/4})^\alpha \le \omega_L
\end{equation*}
and by a simple calculation
$$
(p+q) \left(\theta^N \right)^2 \leq [\theta^{M/4}]^{(\log
[\theta^{M/4}])^\alpha }\le L^{\omega_L}
$$
provided $N$ is large enough. Applying condition
${\bf A_{\boldsymbol \omega}}$ we get that for
all solutions of the equation
\begin{equation} \label{equpq}
\pm j_1 n_{k_1} \pm \dots \pm j_p n_{k_p} \pm j_{p+1} n_{k_{p+1}}
\pm \dots \pm j_{p+q} n_{k_{p+q}} = 0
\end{equation}
with different indices $k_1, \ldots, k_{p+q}$, where
\begin{equation} \label{apq}
1\le j_i \leq (p+q) \left(\theta^N\right)^2,
\end{equation}
the $n_{k_j}$ belong to the $[\theta^{M/4}]$ smallest elements of
$(n_k)_{k \geq 1}$. By construction not a single one of these
elements is contained in $\Delta_M$ or $\Delta_N$. Thus the
equation (\ref{equpq}) subject to (\ref{apq}) has no solution
$\left(k_1, \dots, k_{p+q}\right)$ where $k_1, \dots, k_{p+q}$ are
different and satisfy
$$
k_1,\dots,k_p \in \Delta_M,~k_{p+1},\dots,k_{p+q} \in \Delta_N.
$$
Now
\begin{eqnarray}\label{summ}
& & \E Z_M^p Z_N^q  \nonumber\\
& = & \frac{2^{-p-q}}{\sigma_M^p \sigma_N^q} \sum_{\substack{1
\leq j_q, \dots,j_p \leq [\theta^M]^2,\\1 \leq j_{p+1}
,\dots,j_{p+q} \leq [\theta^N]^2}}
\nonumber\\
& & \qquad \quad \sum_{\substack{k_1, \dots,k_p \in \Delta_M, \\
k_{p+1},\dots,k_{p+q} \in \Delta_N}} a_{j_1} \dots a_{j_{p+q}}
~\mathbf{1} \{\pm j_1 n_{k_1} \pm \dots \pm j_{p+q} n_{k_{p+q}}=0
\}.
\end{eqnarray}
If for some $k_1, \ldots, k_{p+q}$ we have (note that in (\ref{summ})
these indices need not be different)
\begin{equation} \label{pm}
\pm j_1 n_{k_1} \pm \dots \pm j_{p+q} n_{k_{p+q}}=0,
\end{equation}
then grouping the terms of the equation according to identical
indices, we get a new equation of the form
$$
j_1'n_{l_1}+\ldots +j_s' n_{l_s}=0, \qquad l_1<\ldots <l_s, \
s\le p+q, \ \,  j_i'\le (p+q) (\theta^N)^2
$$
and using the above observation, all the coefficients $j_1',
\ldots, j_s'$ must be equal to 0. In other words, in any solution
of (\ref{pm}) the terms can be divided into groups such that in
each group the $n_{k_j}$ are equal and the sum of the coefficients
is 0. Consider first the solutions where all groups have
cardinality 2. This can happen only if both $p$ and $q$ are even,
and similarly to the proof of Lemma \ref{lemma1}, the
contribution of such solutions in (\ref{summ}) is
$$
\frac{p!}{(p/2)!2^{p/2}} \frac{q!}{(q/2)!2^{q/2}}.
$$
Consider now the solutions of (\ref{pm}) where at least one group
has cardinality $\ge 3$. Clearly the sets $\{k_1, \ldots, k_p\}$
and $\{k_{p+1}, \ldots k_{p+q}\}$ are disjoint; let us denote the
number of groups within these two sets by $R$ and $S$,
respectively. Evidently $R\le p/2$, $S\le q/2$, and at least one
of the inequalities is strict. Fixing $j_1, \ldots j_{p+q}$ and the
groups, the number of such solutions cannot exceed
$$|\Delta_M|^R |\Delta_N|^S \le |\Delta_M|^{p/2} |\Delta_N|^{q/2}
 |\Delta_M|^{-1/2}\ll \sigma_M^p\sigma_N^q |\Delta_M|^{-1/2},$$
where we used (\ref{sigmabound}) and the fact that $|\Delta_M|\le |\Delta_N|$.
Since the
number of partitions of the set $\{1, 2, \ldots, p+q\}$ into
disjoint subsets is at most $(p+q)!2^{p+q}$ and since the number
of choices for the signs $\pm$ in (\ref{pm}) is at most $2^{p+q}$,
we see, after summing over all possible values of $j_1, \dots,
j_{p+q}$, that the contribution of the solutions containing at
least one group with cardinality $\ge 3$ in (\ref{summ}) is at most
$2^{p+q} (p+q)! |\Delta_M|^{-1/2} (\log [\theta^{N}])^{p+q}$.
This completes the proof of Lemma~\ref{lemmalil}.

The rest of the proof of the LIL part of Theorem \ref{th1} can be
modeled following the lines of Lemma 4, Lemma 5, Lemma 6 and the
proof of Theorem 1 in \cite{HLP}.

\subsection{Proof of Theorem \ref{th3}}

Let $\omega_k \uparrow \infty$ and
set $\eta_k=\frac{1}{2}\omega_k^{1/2}$,
$\boldsymbol{\eta} = (\eta_1, \eta_2,\dots)$. Clearly
\begin{equation} \label{etak}
(2k)^{\eta_k^2+2 \eta_k} \leq (2k)^{\omega_k/2}
\le k^{-2}|I_k| \qquad \text{for} \ k\ge k_0
\end{equation}
since, as we noted, $|I_k|\ge a \omega_k (k-1)^{\omega_k-1}\ge
(k/2)^{\omega_k-1}$ for large $k$. We choose $n_k$, $k=1, 2,
\ldots$ independently and uniformly from the integers of the
intervals $I_k$ in (\ref{I}). We claim that, with probability 1,
the sequence $(n_k)_{k\ge 1}$ is increasing and satisfies
condition $\boldsymbol{A_\eta}$. To see this, let $k\ge 1$ and
consider the numbers of the form
\begin{equation}\label{const}
(a_1 n_{i_1}+\ldots +a_ s n_{i_s})/d
\end{equation}
where $1 \leq s \leq \eta_k$,  $1\le i_1, \ldots, i_s\le k-1$,
$a_1, \ldots a_s, d$ are nonzero integers with
$|a_1|, \ldots, |a_s|, |d|\le k^{\eta_k}$. Since the number of values in
(\ref{const}) is at most $(2k)^{\eta_k^2+2 \eta_k}$, and by
(\ref{etak}), the probability that $n_k$ equals any of these
numbers is at most $k^{-2}$. Thus by the Borel-Cantelli lemma,
with probability 1 for $k\ge k_1$, $n_k$ will be different from
all the numbers in (\ref{const}) and thus the equation
$$a_1 n_{i_1}+\ldots +a_s n_{i_{s}}+ a_{s+1} n_k=0$$
has no solution with $1 \leq s \leq \eta_k$, $1\le i_1< \ldots <
i_s\le k-1$, $0< |a_1|, \ldots, |a_{s+1}|\le k^{\eta_k}$. By monotonicity,
the equation
$$
a_1n_{i_1}+\ldots + a_s n_{i_s}=0
$$
has no solutions provided the indices $i_\nu$ are all different,
the maximal index is at least $k$,  the number of terms is at most
$\eta_k$, and $0< |a_1|, \ldots, |a_s|\le k^{\eta_k}$. In other
words, $(n_k)$ satisfies condition $\boldsymbol{A_\eta}$.
Using now Theorem \ref{th1}, we get Theorem \ref{th3}.


\end{document}